\documentclass[12pt]{article}

\usepackage{amssymb,amsfonts,amsmath}

\newtheorem{theorem}{Theorem}[section]

\newtheorem{corollary}{Corollary}[section]

\begin{document}

\begin{center}

{\large \bf FUNCTIONAL SPACES AND OPERATORS CONNECTED WITH SOME L\'EVY NOISES}\\[10mm]
\large E. LYTVYNOV\\[2mm]
{\it Department of Mathematics\\
University of Wales Swansea \\
Singleton Park\\
Swansea SA2 8PP\\
U.K.\\
E-mail: e.lytvynov@swansea.ac.uk}

\end{center}

\begin{abstract}

We review some recent developments in white noise
analysis and quantum probability. We pay a special attention to
spaces of test and generalized functionals of some L\'evy white
noises, as well as as to the structure of quantum white noise on
these spaces. \end{abstract}


\section{Gaussian white noise and Fock space}

Since the work of Hida  \cite{hida} of 1975, Gaussian white noise
analysis has become an established theory of test and generalized
functions of infinitely many variables, see e.g.\ 
   \cite{HKPS,BK} and the references therein.

Let us shortly recall some basic results of Gaussian analysis. In
the space $L^2({\mathbb R}):=L^2({\mathbb R},dx)$, consider the
harmonic oscillator
$$ (Hf)(t):=-f''(t)+(t^2+1)f(t),\qquad f\in C_0^\infty({\mathbb R}).$$ This operator is self-adjoint and we preserve the notation $H$ for its closure. For each $p\in{\mathbb R}$, define a scalar product
$$ (f,g)_p:=(H^pf,g)_{L^2({\mathbb R})},\qquad f,g\in C_0^\infty({\mathbb R}).$$
Let $\mathcal S_p$ denote the Hilbert space obtained as the closure
of $C_0^\infty({\mathbb R})$ in the norm $\|\cdot\|_p$ generated by
the scalar product $(\cdot,\cdot)_p$. Then, for any $p>q$, the space
$\mathcal S_p$ is densely and continuously embedded into $\mathcal
S_q$, and if $p-q>1/2$, then this embedding is  of Hilbert--Schmidt
type. Furthermore, for each $p>0$, $\mathcal S_{-p}$ is the dual
space of $\mathcal S_p$ with respect to zero space $L^2({\mathbb
R})$, i.e., the dual paring $\langle f,\varphi\rangle$ between any
$f\in \mathcal S_{-p}$ and any $\varphi\in \mathcal S_{p}$ is
obtained as the extension of the scalar product in $L^2({\mathbb
R})$. The above conclusions are, in fact, corollaries of the fact
that the sequence of Hermite functions on ${\mathbb R}$,
$$ e_j=e_j(t)=(\sqrt\pi 2^jj!)^{-1/2}(-1)^j e^{t^2/2}(d/dt)^je^{-t^2},\qquad j\in{\mathbb Z}_+:=\{0,1,2,\dots\},$$
forms an orthonormal basis of $L^2({\mathbb R})$ such that each
$e_j$ is an eigenvector of $H$ with eigenvalue $(2j+2)^2$.

Then $$\mathcal S := \projlim _{p\to\infty}\mathcal S_p$$ is the
Schwartz space of infinitely differentiable, rapidly decreasing
functions on ${\mathbb R}$, and its dual $$ \mathcal
S'=\operatornamewithlimits{ind\,lim}_{p\to\infty} \mathcal S_{-p} $$
is the Schwartz space of tempered distributions.

We denote by $\mathcal C({\mathcal S'})$ the $\sigma$-algebra on
$\mathcal S'$ which is generated by cylinder sets in $\mathcal S'$,
i.e., by the sets of the form $$\{\omega\in\mathcal S':
(\langle\omega_,\varphi_1\rangle,\dots,\langle
\omega,\varphi_N\rangle)\in A\},$$ where
$\varphi_1,\dots,\varphi_N\in\mathcal S$, $N\in{\mathbb N}$, and
$A\in {\mathcal B}({\mathbb R}^N)$.

By the Minlos theorem, there exists a unique probability measure
$\mu_{\mathrm G}$ on $(\mathcal S',\mathcal C(\mathcal S'))$ whose
Fourier transform is given by \begin{equation}\label{jhf}
\int_{\mathcal S'}e^{i\langle\omega,\varphi\rangle}d\mu_{\mathrm
G}(\omega) =\exp\big[-(1/2)\|\varphi\|_{L^2({\mathbb
R})}^2\big],\qquad \varphi\in\mathcal S.\end{equation}

The measure $\mu_{\mathrm G}$ is called the (Gaussian) white noise
measure. Indeed, using formula \eqref{jhf}, it is easy to see that,
for each $\varphi\in\mathcal S$, $$ \int_{\mathcal
S'}\langle\omega,\varphi\rangle^2\,d\mu_{\mathrm
G}(\omega)=\|\varphi\|^2_{L^2({\mathbb R})}.$$ Hence, extending the
mapping
$$ L^2({\mathbb R})\supset\mathcal S \ni\varphi\mapsto \langle\cdot,\varphi\rangle\in
L^2(\mathcal S',\mu_{\mathrm G})$$ by continuity, we obtain a random
variable $\langle\cdot,f\rangle\in L^2(\mathcal S',\mu_{\mathrm G})$
for each $f\in L^2({\mathbb R})$. Then, for each $t\in{\mathbb R}$,
we define
\begin{equation}\label{ifulfr} X_t:=\begin{cases}
\langle\cdot,\pmb1_{[0,t]}\rangle,&t\ge0,\\
-\langle\cdot,\pmb1_{[t,0]}\rangle,&t<0.\end{cases}
\end{equation} It is easily seen that $(X_t)_{t\in{\mathbb R}}$ is a version of Brownian motion, i.e., finite-dimensional distributions of the stochastic process $(X_t)_{t\in{\mathbb R}}$ coincide with those of Brownian motion. We now
informally have, for all $t\in{\mathbb R}$,
$X_t(\omega)=\int_0^t\omega(t)\,dt$, so that
$X_t'(\omega)=\omega(t)$. Thus, elements $\omega\in\mathcal S'$ can
be thought of as  paths of the derivative of Brownian motion, i.e.,
Gaussian white noise.

Let us recall that the symmetric  Fock space over a real separable
Hilbert space $\mathcal H$ is defined as
$$\mathcal F(\mathcal H):=\bigoplus_{n=0}^\infty
{\mathcal F}^{(n)}(\mathcal H)n!\, .$$ Here  ${\mathcal
F}^{(n)}(\mathcal H):=\mathcal H^{\odot n}_{\mathbb C}$, where
$\odot$ stands for symmetric tensor product and the lower index
${\mathbb C}$ denotes complexification of a real space. Thus, for
each $(f^{(n)})_{n=0}^\infty\in\mathcal F(\mathcal H)$, $$
\|(f^{(n)})_{n=0}^\infty\|_{\mathcal F(\mathcal H)}^2=
\sum_{n=0}^\infty \|f^{(n)}\|_{\mathcal F^{(n)}(\mathcal H)}^2n!\, .
$$

The central technical point of the construction of spaces of test
and generalized  functionals of Gaussian white noise is the
Wiener--It\^o--Segal isomorphism $I_{\mathrm G}$ between the Fock
space $\mathcal F(L^2({\mathbb R}))$ and the complex space
$L^2(\mathcal S'\to{\mathbb C},\mu) $, which, for simplicity of
notations, we will denote by $L^2(\mathcal S',\mu) $.

There are different ways of construction of  the isomorphism
$I_{\mathrm G}$, e.g., using multiple stochastic integrals with
respect to Gaussian random measure. For us, it will be convenient to
follow the approach which  uses the procedure of orthogonalization
of polynomials, see e.g.\  \cite{BK} for details.

A function $F(\omega)=\sum_{i=0}^n\langle\omega^{\otimes
i},f^{(i)}\rangle$, where $\omega\in\mathcal S'$, $n\in{\mathbb
Z}_+$, and each $f^{(i)}\in\mathcal S^{\odot i}_{\mathbb C}$, is
called a continuous polynomial on $\mathcal S'$, and $n$ is called
the order of the polynomial $F$. The set $\mathcal P$ of all
continuous polynomials on $\mathcal S'$ is dense in $L^2({\mathcal
S}',\mu_{\mathrm G})$. For $n\in{\mathbb Z}_+$, let $\mathcal
P^{(n)}$ denote the set of all continuous polynomials on $\mathcal
S'$ of order $\le n$, and let $\mathcal P^{(n)}_{\mathrm G}$ be the
closure of $\mathcal P^{(n)}$ in $L^2(\mathcal S',\mu_{\mathrm G})$.
Let $\mathfrak P_{{\mathrm G}}^{(n)}$ stand for the orthogonal
difference $\mathcal P^{(n)}_{{\mathrm G}}\ominus \mathcal
P^{(n-1)}_{{\mathrm G}}$ in $L^2(\mathcal S',\mu_{\mathrm G})$. Then
we easily get the orthogonal decomposition $$ L^2(\mathcal
S',\mu_{\mathrm G})=\bigoplus_{n=0}^\infty \mathfrak
P^{(n)}_{{\mathrm G}}. $$

Next, for any $f^{(n)}\in \mathcal S^{\odot n}_{\mathbb C}$, we
define ${:}\langle \omega^{\otimes n},f^{(n)}\rangle{:}_{{\mathrm
G}}$ as the orthogonal projection of $\langle \omega^{\otimes
n},f^{(n)}\rangle$ onto $\mathfrak P^{(n)}_{{\mathrm G}}$. The set
of such projections is dense in $\mathfrak P^{(n)}_{{\mathrm G}}$.
Furthermore, for any $f^{(n)},g^{(n)}\in\mathcal S^{\odot
n}_{\mathbb C}$, we have:
\begin{equation}\label{hhgcv} \int_{\mathcal S'}\overline{{:}\langle
\omega^{\otimes n},f^{(n)}\rangle{:}_{{\mathrm G}}}\times {:}\langle
\omega^{\otimes n},g^{(n)}\rangle{:}_{{\mathrm G}}\,d\mu_{\mathrm
G}(\omega)=(f^{(n)},g^{(n)})_{\mathcal F^{(n)}(L^2({\mathbb
R}))}n!\,.
\end{equation}

Let $\mathcal F_{\mathrm{fin}}(\mathcal S)$ denote the set of all
sequences $(f^{(n)})_{n=0}^\infty$  such that each
$f^{(i)}\in\mathcal S^{\odot i}_{\mathbb C}$ and for some
$N\in{\mathbb N}$ $f^{(n)}=0$ for all $n\ge N$.
 The  $\mathcal F_{\mathrm{fin}}(\mathcal S)$ is a dense subset of $\mathcal F(L^2({\mathbb R}))$. For any $f=(f^{(n)})_{n=0}^\infty\in\mathcal
  F_{\mathrm{fin}}(\mathcal S)$, we set \begin{equation}\label{jhfhfc}
(I_{{\mathrm G}}f)(\omega)=\sum_{n=0}^\infty
{:}\langle\omega^{\otimes n}, f^{(n)}\rangle{:}_{{\mathrm G}}\in
L^2(\mathcal S',\mu_{\mathrm G}) .\end{equation} By \eqref{hhgcv},
we can extend $I_{{\mathrm G}}$ by continuity to get a unitary
operator
$$I_{{\mathrm G}}:\mathcal F(L^2({\mathbb R}))\to L^2(\mathcal S',\mu_{\mathrm
G} ).$$  For any function $f^{(n)}\in\mathcal F^{(n)}(L^2({\mathbb
R}))$, we will use the evident notation ${:}\langle\omega^{\otimes
n}, f^{(n)}\rangle{:}_{{\mathrm G}}$. Then, for each
$f=(f^{(n)})_{n=0}^\infty\in \mathcal F(L^2({\mathbb R}))$,
$I_{{\mathrm G}}f$ is given by  formula \eqref{jhfhfc}.

For any  $\varphi\in\mathcal S$, let
$\langle\cdot,\varphi\rangle\cdot$ denote the operator of
multiplication by $\langle\cdot,\varphi\rangle$ in $L^2(\mathcal
S',\mu_{\mathrm G})$. We set $$ A_{{\mathrm G}}(\varphi):=
I_{\mathrm G}^{-1}\langle\cdot,\varphi\rangle\cdot I_{\mathrm G}.$$
Then, $\mathcal F_{\mathrm{fin}}(\mathcal S)\subset
\operatorname{Dom}(A_{{\mathrm G}}(\varphi))$, $A_{{\mathrm
G}}(\varphi)\mathcal F_{\mathrm{fin}}(\mathcal S)\subset \mathcal
F_{\mathrm{fin}}(\mathcal S)$, and $A_{{\mathrm G}}(\varphi)$ is
essentially self-adjoint on $\mathcal F_{\mathrm{fin}}(\mathcal S)$.
Furthermore, $A_{{\mathrm G}}(\varphi)$ has the following
representation on $\mathcal F_{\mathrm{fin}}(\mathcal S)$:
\begin{equation}A_{{\mathrm G}}(\varphi)=A^+(\varphi)+A^-(\varphi).
\label{fgffg}\end{equation} Here $A^+(\varphi)$ is the creation
operator: for $f^{(n)}\in\mathcal S^{\odot n}_{\mathbb C}$
$$A^+(\varphi)f^{(n)}=\varphi\odot f^{(n)}\in \mathcal
S^{\odot(n+1)}_{\mathbb C},$$ and $A^-(\varphi)$ is the annihilation
operator: $A^-(\varphi)f^{(n)}$ belongs to $\mathcal S_{\mathbb
C}^{\odot(n-1)}$ and is given by
$$ (A^-(\varphi)f^{(n)})(t_1,\dots,t_{n-1})=\int_{\mathbb R} \varphi(t)f^{(n)}(t,t_1,\dots,t_{n-1})\,dt.$$

For each $\varkappa\in[-1,1]$ and $p\in{\mathbb R}$, we denote
$$
\mathcal F_\varkappa(\mathcal S_p):=\bigoplus_{n=0}^\infty \mathcal
F^{(n)}(\mathcal S_p)(n!)^{1+\varkappa},
$$
and for each $\varkappa\in[0,1]$ we set $$ \mathcal F_\varkappa
(\mathcal S):=\projlim_{p\to\infty} \mathcal F_{\varkappa}(\mathcal
S_p).$$ It is easy to see that each $\mathcal F_\varkappa (\mathcal
S)$ is a nuclear space. Furthermore, the dual space of $\mathcal
F_\varkappa (\mathcal S)$ with respect to zero space $\mathcal
F(L^2({\mathbb R}))$ is $$ \mathcal F_{-\varkappa}(\mathcal
S'):=\operatornamewithlimits{ind\,lim}_{p\to\infty }\mathcal F
_{-\varkappa}(\mathcal S_{-p}).$$ Thus, we get a standard triple
$$ \mathcal F_{\varkappa}(\mathcal S)\subset \mathcal F
(L^2({\mathbb R}))\subset \mathcal F_{-\varkappa}(\mathcal S').$$
The test space $\mathcal F_1(\mathcal S)$ is evidently the smallest
one between the above spaces, whereas its dual space $\mathcal
F_{-1}(\mathcal S')$ is the biggest one.

For each $F=(F^{(n)})\in\mathcal F_{-1}(\mathcal S')$, the
$S$-transform of $F$ is defined by \begin{equation}\label{dters}
(SF)(\xi):= \sum_{n=0}^\infty \langle
\overline{F^{(n)}},\xi^{\otimes n}\rangle,\qquad \xi\in\mathcal
S_{\mathbb C},\end{equation} provided the series on the right hand
side of \eqref{dters} converges absolutely.

The $S$-transform uniquely characterizes an element of $\mathcal
F_{-1}(\mathcal S')$.
 More exactly, let $\operatorname{Hol}_0(\mathcal S_{{\mathbb C}})$ denotes the set of all (germs) of functions which are holomorphic in a neighborhood of zero in $\mathcal S_{{\mathbb C}}$.
The following theorem was proved in \cite{KLS}.

\begin{theorem}\label{hjgv} The $S$-transform is a one-to-one map between  $\mathcal F_{-1}(\mathcal S')$ and  $\operatorname{Hol}_0(\mathcal S_{{\mathbb C}})$.

\end{theorem}

Note that the choice of $\varkappa>1$ would imply that the
$S$-transform is not well-defined on $\mathcal
F_{-\varkappa}(\mathcal S')$.  On the other hand, all the spaces
$\mathcal F_{-\varkappa}(\mathcal S')$ and $\mathcal
F_{\varkappa}(\mathcal S)$ with $\varkappa\in[0,1]$ admit a complete
characterization in terms of their $S$-transform, see e.g.\    \cite{BK,KS,KLS}.

Taking into account that the product of two elements of
$\operatorname{Hol}_0(\mathcal S_{{\mathbb C}})$ remains in this
set, one defines a Wick product $F_1\diamond F_2$ of $F_1,F_2\in
\mathcal F_{-1}(\mathcal S')$ through the formula
\begin{equation}\label{aewra}S(F_1\diamond
F_2)=S(F_1)S(F_2).\end{equation} Furthermore, if $F\in \mathcal
F_{-1}(\mathcal S')$ and $f$ is a holomorphic function in a
neighborhood of $(SF)(0)$ in ${\mathbb C}$, then one defines
$f^\diamond(F)\in\mathcal F_{-1}(\mathcal S')$ through
\begin{equation}\label{ufu} S(f^\diamond(F))=f(SF).\end{equation}

Using the unitary operator $I_{{\mathrm G}}$, all the above
definitions and results can be reformulated in terms of test and
generalized functions on $\mathcal S'$ whose dual paring is
generated by the scalar product in $L^2(\mathcal S',\mu_{\mathrm
G})$. In particular, one defines spaces of test functions $
(\mathcal S)_{{\mathrm G}}^\varkappa:=I_{{\mathrm G}}\mathcal
F_{\varkappa}(\mathcal S)$ and their dual spaces $(\mathcal
S)_{{\mathrm G}}^{-\varkappa}$, $\varkappa\in[0,1]$. For
$\varkappa=0$, these are the Hida test space and the space of Hida
distrubutions, respectively (e.g.\  \cite{HKPS,BK}).  For
$\varkappa\in(0,1)$, these spaces were introduced and studied by
Kondratiev and Streit   \cite{KS}, and for $\varkappa=1$, by
Kondratiev, Leukert and Streit  \cite{KLS}. Note also that, in the
case of a Gaussian product measure, such spaces for all
$\varkappa\in[0,1]$  were studied by Kondratiev  \cite{K}.

The Wick calculus of generalized Gaussian functionals based on the
definitions  \eqref{aewra}, \eqref{ufu}  and the unitary operator
$I_{\mathrm G}$ has found numerous applications, in particular, in
fluid mechanics and financial mathematics, see e.g.\
  \cite{Oks1,EV}.

Additionally to the description of the above test spaces $(\mathcal
S)_{\mathrm G}^\varkappa$ in terms of their $S$-transform, one can
also give their inner description, e.g.\     \cite{KS,KLS}. Let
$\mathcal E(\mathcal S_{{\mathbb C}}')$ denote the space of all
entire functions on $\mathcal S'_{\mathbb C}$. For each
$\beta\in[1,2]$, we denote by $\mathcal
E_{\mathrm{min}}^\beta(\mathcal S_{{\mathbb C}}')$ the subset of
$\mathcal E(\mathcal S_{{\mathbb C}}')$ consisting of all entire
functions of the $\beta$-th order of growth and minimal type. That
is, for any $\Phi\in \mathcal E_{\mathrm{min}}^\beta(\mathcal
S_{{\mathbb C}}')$, $p\ge0$, and $\varepsilon>0$, there exists $C>0$
such that
$$|\Phi(z)|\le C\exp(\varepsilon |z|_{-p}^\beta),\qquad z\in\mathcal S_{-p,{\mathbb C}}.$$ Next, we denote by $\mathcal E_{\mathrm{min}}^\beta(\mathcal S')$ the set of all functions on $\mathcal S'$, which are obtained by restricting functions from $\mathcal E_{\mathrm{min}}^\beta(\mathcal S_{{\mathbb C}}')$ to $\mathcal S'$.
The following theorem unifies the results of  \cite{KS,KLS},
see also     \cite{BK,HKPS},

\begin{theorem}\label{theorem}
For each $\varkappa\in[0,1]$, we have
\begin{equation}\label{tdrt} (\mathcal S)_{\mathrm
G}^\varkappa=\mathcal E_{\mathrm{min}}^{2/(1+\varkappa)}(\mathcal
S').\end{equation} The equality \eqref{tdrt} is understood in the
sense that, for each $f=(f^{(n)})\in\mathcal F_{\varkappa}(\mathcal
S)$, the following realization was chosen for $\Phi:=I_{\mathrm
G}f$:
$$\Phi(\omega)= \sum_{n=0}^\infty \langle{:}\omega^{\otimes n}{:}_{\mathrm G},f^{(n)}\rangle,\qquad\omega\in\mathcal S',$$
where ${:}\omega^{\otimes n}{:}_{\mathrm G}\in\mathcal S'{}^{\odot
n}$ is defined by the recurrence relation \begin{gather*}
{:}\omega^{\otimes 0}{:}_{\mathrm G}=1,\quad  {:}\omega^{\otimes 1}{:}_{\mathrm G}=\omega,\\
{:}\omega^{\otimes (n+1)}{:}_{\mathrm G}(t_1,\dots,t_{n+1})=\big(
{:}\omega^{\otimes n}{:}_{\mathrm G}(t_1,\dots,t_n)\omega(t_{n+1})\big)^\sim \\
\text{}-n \big({:}\omega^{\otimes(n-1)}{:}_{\mathrm
G}(t_1,\dots,t_{n-1})1(t_n)\delta(t_{n+1}-t_{n})\big)^\sim,\qquad
n\in{\mathbb N}. \end{gather*}Here $\delta(\cdot)$ denotes the delta
function at zero and $(\cdot)^\sim$ denotes symmetrization.
\end{theorem}

For each $t\in{\mathbb R}$, we define an annihilation operator at $t$, denoted by $\partial _t$, and  a creation operator at $t$, denoted by $\partial ^\dag_t$,   by \begin{align*} (\partial_t f^{(n)} )(t_1,\dots,t_{n-1})&:=nf^{(n)}(t_1,\dots,t_{n-1},t),\qquad f^{(n)}\in\mathcal S_{\mathbb C}^{\odot n}, \\
\partial_t^\dag F^{(n)}&:=\delta_t\odot F^{(n)},\qquad F^{(n)}\in\mathcal
S_{\mathbb C}^{\prime\,\odot n},\end{align*} where $\delta_t$
denotes the delta function at $t$. The operators $\partial_t$ and
$\partial_t^\dag$ can then be  extended to linear continuous
operators on $\mathcal F_{\varkappa}(\mathcal S)$ and $\mathcal
F_{-\varkappa}(\mathcal S')$, respectively, and $\partial ^\dag_t$
becomes the dual operator of $\partial_t$. Then,
\begin{align*}
A^-(\varphi)&=\int_{\mathbb R}\varphi(t)\partial_t\,dt,\\
A^+(\varphi)&=\int_{\mathbb R}\varphi(t)\partial_t^\dag\,dt,
\end{align*}
and so   by \eqref{fgffg}, $$ A_{\mathrm G}(\varphi)=\int_{\mathbb
R} \varphi(t)(\partial_t+\partial_t^\dag)\,dt,$$ the above integrals
being understood in the Bochner sense. The operators
\begin{equation}\label{jft} W_{\mathrm
G}(t):=\partial_t+\partial_t^\dag,\qquad t\in{\mathbb
R},\end{equation} are called quantum Gaussian white noise. Realized
on $(\mathcal S)_{\mathrm G}^1$, the operator $\partial_t$ becomes
the operator of G\^ateaux differentiation in direction $\delta_t$:
$$ (\partial_t F)(\omega)=\lim_{\varepsilon\to0}
(F(\omega+\varepsilon\delta_t)-F(\omega))/\varepsilon,\qquad F\in
(\mathcal S)_{\mathrm G}^1,\ t\in{\mathbb R},\ \omega\in\mathcal S',
$$ see e.g.\   \cite{HKPS}.

\section{Poisson white noise} The measure $\mu_{\mathrm P}$ of (centered) Poisson white noise is defined on $(\mathcal S',\mathcal C(\mathcal S'))$ by $$
\int_{\mathcal S'}e^{i\langle \omega,\varphi\rangle}\,d\mu_{\mathrm
P}(\omega) =\exp\bigg[\int_{\mathbb
R}(e^{i\varphi(t)}-1-i\varphi(t))\,dt\bigg],\qquad
\varphi\in\mathcal S.$$

Under $\mu_{\mathrm P}$, the stochastic process $(X_t)_{t\in{\mathbb
R}}$, defined by \eqref{ifulfr}, is a centered Poisson process,
which is why $\omega\in\mathcal S'$ can now be thought of as a path
of Poisson white noise (in fact, $\mu_{\mathrm P}$ is concentrated
on infinite sums of delta functions shifted by $-1$).

The procedure of orthogonalization of continuous polynomials in
$L^2(\mathcal S',\mu_{\mathrm P})$ leads to a unitary isomorphism
$I_{\mathrm P}$ between the Fock space $\mathcal F(L^2({\mathbb
R}))$ and $L^2(\mathcal S',\mu_{\mathrm P})$. The counterpart of
formula \eqref{fgffg} now looks as follows:
$$A_{{\mathrm P}}(\varphi)=A^+(\varphi)+A^0(\varphi)+A^-(\varphi),$$
where $A^0(\varphi)$ is the neutral operator:
$$ (A^0(\varphi)f^{(n)}) (t_1,\dots,t_n)
:=\bigg(\sum_{i=1}^n \varphi(t_i)\bigg)f^{(n)}(t_1,\dots,t_n),\qquad
f^{(n)} \in\mathcal S_{\mathbb C}^{\odot n.}$$ In terms of the
operators $\partial_t$ and $\partial_t^\dag$, the neutral operator
has the following representation:
$$ A^0(\varphi)=\int_{\mathbb R}\varphi(t)\partial_t^\dag\partial_t\,dt.$$
Thus, the quantum (centered) Poisson white noise is given by
\begin{equation}\label{jgufgujjj}W_{\mathrm P}(t)=\partial_t+\partial_t^\dag\partial_t+\partial_t^\dag\end{equation}
(see e.g.\     \cite{IK,L2} ).

Next, using the unitary operator $I_{\mathrm P}$, we obtain a scale
of spaces of test functions $(\mathcal S)_{\mathrm P}^\varkappa$ and
generalized functions $(\mathcal S')_{\mathrm P}^{-\varkappa}$. For
any $f^{(n)}\in\mathcal S_{\mathbb C}^{\odot n}$, the orthogonal
projection ${:}\langle\omega^{\otimes n},f^{(n)}\rangle{:}_{\mathrm
P}$ has a $\mu_{\mathrm P}$-version $\langle {:}\omega^{\otimes
n}{:}_{\mathrm P},f^{(n)}\rangle$, where ${:}\omega^{\otimes
n}{:}_{\mathrm P}\in\mathcal S'{}^{\odot n}$ are given by the
following recurrence relation (see   \cite{L1}):
\begin{gather*}
{:}\omega^{\otimes 0}{:}_{\mathrm P}=1,\quad  {:}\omega^{\otimes 1}{:}_{\mathrm P}=\omega,\\
{:}\omega^{\otimes(n+1)}{:}_{\mathrm P}(t_1,\dots,t_{n+1})=\big(
{:}\omega^{\otimes n}{:}_{\mathrm P}(t_1,\dots,t_n)\omega(t_{n+1})\big)^\sim \\
\text{}-n \big({:}\omega^{\otimes(n-1)}{:}_{\mathrm P}(t_1,\dots,t_{n-1})1(t_n)\delta(t_{n+1}-t_{n})\big)^\sim\\
\text{}-n \big({:}\omega^{\otimes n}{:}_{\mathrm
P}(t_1,\dots,t_{n})\delta(t_{n+1}-t_{n})\big)^\sim,\qquad
n\in{\mathbb N}.
\end{gather*}

However, in the Poisson case, the following statement
holds  \cite{L1}.

\begin{theorem}\label{ufytd} For each $\omega\in\mathcal S'$, denote $D(\omega):=({:}\omega^{\otimes n}{:}_{\mathrm P})_{n=0}^\infty$. Then
$D(\omega)\in\mathcal F_{-1}(\mathcal S')$, and if $\omega\ne0$,
then  $D(\omega)\not\in\mathcal F_{-\varkappa}(\mathcal S')$ for all
$\varkappa\in[0,1)$.
\end{theorem}

It is straightforward to see that the $D(\omega)$ in the Poisson
realization is just the delta function at $\omega$, denoted by
$\delta_\omega$. Thus, Theorem~\ref{ufytd} implies:

\begin{corollary}\label{ytdster}
For each $\omega\in\mathcal S'$, $\delta_\omega\in(\mathcal
S')_{\mathrm P}^{-1}$, and if $\omega\ne0$, then
$\delta_\omega\not\in(\mathcal S')_{\mathrm P}^{-\varkappa}$ for all
$\varkappa\in[0,1)$.
\end{corollary}

The above corollary shows that the test spaces $(\mathcal
S)_{\mathrm P}^\varkappa$ with $\varkappa\in[0,1)$ do not posses
nice inner properties, and therefore they are not appropiate for
applications. On the other hand, we have  \cite{KSWY}:

\begin{theorem}\label{hyguy}
We have:
$$ (\mathcal S)_{\mathrm P}^1=\mathcal E_{\mathrm{min}}^1(\mathcal S').$$
\end{theorem}

Thus, $\mathcal E_{\mathrm{min}}^1(\mathcal S')$ appears to be a
universal space for both Gaussian and Poisson white noise analysis.

The annihilation operator $\partial_t$ realized on $(\mathcal
S)_{\mathrm P}^1$ becomes a difference operator  \cite{L1}:
$$ (\partial_t F)(\omega)=F(\omega+\delta_t)-F(\omega),\qquad F\in
(\mathcal S)_{\mathrm P}^1,\ t\in{\mathbb R},\ \omega\in\mathcal
S'.$$

We note that one can also study a more general white noise measure
of Poisson type for which the corresponding quantum white noise is
given by
\begin{equation}\label{gyul}\partial_t+\lambda\partial_t^\dag\partial_t+\partial_t^\dag,\end{equation}
where $\lambda\in{\mathbb R}_+$, $\lambda\ne0$, see e.g.\
 \cite{LRS}.

\section{L\'evy white noise and extended Fock space}\label{hgyd}

We will now discuss the case of a L\'evy white noise without
Gaussian part. Let $\nu$ be a probability measure on $({\mathbb
R},\mathcal B({\mathbb R}))$ such that $\nu(\{0\})=0$. We will also
assume that there exists $\varepsilon>0$ such that $\int_{\mathbb
R}\exp(\varepsilon|s|)\,\nu(ds)<\infty$. The latter condition
implies that $\nu$ has all moments finite, and moreover, the set of
all polynomials on ${\mathbb R}$ is dense in $L^2({\mathbb R},\nu)$.

We define a centered L\'evy white noise as a probability measure
$\mu_\nu$ on $(\mathcal S',\mathcal C(\mathcal S'))$ with Fourier
transform $$\int_{\mathcal S'}e^{i\langle
\omega,\varphi\rangle}\,\mu_\nu(d\omega)=\exp\bigg[\varkappa
\int_{\mathbb R}\int_{\mathbb R}
(e^{is\varphi(t)}-1-is\varphi(t))\frac1{s^2}\,\nu(ds)\,dt\bigg],\qquad
\varphi\in\mathcal S,$$ where $\varkappa>0$. For notational
simplicity, we will assume that $\varkappa=1$, which is a very weak
restriction.

Under the measure $\mu_\nu$, the stochastic process
$(X_t)_{t\in{\mathbb R}}$, defined by \eqref{ifulfr}, is a centered
L\'evy process with L\'evy measure $\frac1{s^2}\,\nu(ds)$:
$$ \int_{\mathcal S'}e^{iu X_t(\omega)}\,\mu_\nu(d\omega)=\exp\bigg[|t|\int_{\mathbb R} (e^{isu\operatorname{sign}(t)}-1-isu\operatorname{sign}(t))\frac1{s^2}\,\nu(ds)\bigg],\qquad u\in{\mathbb R}.$$
Hence, $\omega\in\mathcal S'$ can be thought of as a path of L\'evy
white noise.

By the above assumptions, the set $\mathcal P$ of all continuous
polynomials on $\mathcal S'$ is dense in $L^2(\mathcal S',\mu_\nu)$.
Therefore, through the procedure of orthogonalization  of
polynomials, one gets an orthogonal decomposition
$$ L^2(\mathcal S',\mu_\nu)=\bigoplus _{n=0}^\infty \mathfrak P_\nu^{(n)},$$
and the set of all orthogonal projections ${:}\langle
\omega^{\otimes n},f^{(n)}\rangle{:}_\nu$ of $\langle
\omega^{\otimes n},f^{(n)}\rangle$ onto $\mathfrak P_\nu^{(n)}$ is
dense in $\mathfrak  P_\nu^{(n)}$.

However, in contrast to the Gaussian and Poisson  cases, the scalar
product of any ${:}\langle \omega^{\otimes n},f^{(n)}\rangle{:}_\nu$
and ${:}\langle \omega^{\otimes n},g^{(n)}\rangle{:}_\nu$ in
$L^2(\mathcal S',\mu_\nu)$ is not given by the scalar product of
 $f^{(n)}$ and $g^{(n)}$ in the Fock space, but by a much more
 complex expression, see     \cite{L4,BLM} for an explicit formula in the general case,
  and formulas \eqref{1111}, \eqref{2222} below in a special case. Still it is possible to construct a unitary isomorphism $I_\nu$ between the so-called extended Fock space $F_\nu(L^2({\mathbb R}))=\bigoplus_{n=0}^\infty F_\nu^{(n)}(L^2({\mathbb R}))n!$ and $L^2(\mathcal S',\mu_\nu)$. In the above construction, $F_\nu^{(n)}(L^2({\mathbb R}))$ is a Hilbert space that is obtained as the closure of $\mathcal S_{\mathbb C}^{\odot n}$ in the norm generated by the scalar product $$ (f^{(n)},g^{(n)})_{F_\nu^{(n)}(L^2({\mathbb R}))}:=\frac1{n!}\int_{\mathcal S'}
\overline{{:}\langle \omega^{\otimes n},f^{(n)}\rangle
{:}_\nu}\times {:}\langle \omega^{\otimes n},g^{(n)}\rangle {:}_\nu
\,\mu_\nu(d\omega). $$

We next set $$ A_\nu(\varphi):=I_\nu^{-1}\langle
\cdot,\varphi\rangle\cdot I_\nu,\qquad \varphi\in\mathcal S.$$ Our
next aim is to derive an explicit form of the action of these
operators. This can be done  \cite{L4,BLM}, however, the property that
$\mathcal F_{\mathrm{fin}}(\mathcal S)$ is invariant under the
action of $A_\nu(\varphi)$, generally speaking, does not hold. The
following theorem  \cite{BLM,L4,L3} identifies all the L\'evy noises
for which this property is preserved.

\begin{theorem}\label{hfydytd}
Under the above assumptions, the property
$$ A_\nu(\varphi)\mathcal F_{\mathrm{fin}}(\mathcal S)\subset \mathcal F_{\mathrm{fin}}(\mathcal S),\qquad \varphi\in\mathcal S,$$
holds if and only if $\nu$ is the measure of orthogonality of a
system of polynomials $(p_n(t))_{n=0}^\infty$ on ${\mathbb R}$ which
satisfy the following recurrence relation:
\begin{gather}
tp_n(t)=\sqrt{(n+1)(n+2)}p_{n+1}(t)+\lambda
(n+1)p_n(t)+\sqrt{n(n+1)} p_{n-1}(t), \notag \\ n\in{\mathbb Z}_+,\
p_0(t)=1,\ p_{-1}(t)=0,\label{uftuy}
\end{gather}
for some $\lambda\in{\mathbb R}$. \end{theorem}

Let us consider the situation described in Theorem~\ref{hfydytd}
 in more detail. We will denote by $\nu_\lambda$ the
  measure $\nu$ which corresponds to the parameter
  $\lambda\in{\mathbb R}$ through \eqref{uftuy}.
  We first mention that the condition of orthogonality of
  the polynomials satisfying \eqref{uftuy} uniquely determines
   the measure $\nu_\lambda$.
   It is easy to show that, for $\lambda>0$,
   the measure $\nu_{-\lambda}$ is the image of the measure
   $\nu_\lambda$ under the mapping ${\mathbb R}\ni t\mapsto -t\in{\mathbb R}$,
   which is why we will only consider the case $\lambda\ge0$.
   In fact, we have (see e.g.\ Ref.~   \cite{Chihara}):
   for $\lambda\in[0,2)$, \begin{align*} \nu_\lambda(ds)&=
\frac{\sqrt{4-\lambda^2}}{2\pi}  \big|\Gamma\big(1+i
(4-\lambda^2)^{-1/2}s\big)\big|^2\\
&\times\exp\big[-s2(4-\lambda^2)^{-1/2}\arctan
\big(\lambda(4-\lambda^2)^{-1/2}\big) \big]\,ds \end{align*}
($\nu_\lambda$ is a Meixner distribution), for $\lambda=2$
$$\nu_2(ds)=\chi_{(0,\infty)} (s)e^{-s}s\,ds$$
($\nu_2$ is a gamma distribution), and for $\lambda>2$
$$\nu_\lambda(ds)=(\lambda^2-4)\sum_{k=1}^\infty
\bigg(\frac{\lambda-\sqrt{\lambda^2-4}}{\lambda+\sqrt{\lambda^2-4}}\bigg)^k\,k\,\delta_{\sqrt{\lambda^2-4}\,k}$$
($\nu_\lambda$ is now a Pascal distribution).

In what follows,we will use the lower index $\lambda$ instead of
$\nu_{\lambda}$.

The stochastic process $(X_t)_{t\in{\mathbb R}}$ under $\mu_\lambda$
is a Meixner process for $|\lambda|<2$, a gamma process for
$|\lambda|=2$ and a Pascal process  for $|\lambda|>2$. In other
words, for each $t\in{\mathbb R}$, $t\ne0$, the distribution of the
random variable $X_t$ under $\mu_{\lambda}$ is of the same class of
distributions as the measure $\nu_\lambda$.

Next, for each $\lambda$, the scalar product in the space
$F^{(n)}(L^2({\mathbb R}))=F^{(n)}_{\lambda}(L^2({\mathbb R}))$ is
given as follows  \cite{L3,L4}. For each $\alpha\in{\mathbb Z}_+$,
$1\alpha_1+2\alpha_2+\dots=n$, $n\in{\mathbb N}$, and for any
function $f^{(n)}:{\mathbb R}^n\to{\mathbb R}$ we define a function
$D_\alpha f^{(n)}:{\mathbb R}^{|\alpha|}\to{\mathbb R}$ by setting
\begin{align}(D_\alpha f^{(n)})(t_1,\dots,t_{|\alpha|}){:=}&
f^{(n)}(t_1,\dots,t_{\alpha_1},
\underbrace{t_{\alpha_1+1},t_{\alpha_1+1}}_{\text{2 times }},
\underbrace{t_{\alpha_1+2},t_{\alpha_1+2}}_{\text{2 times }},\dots,
\underbrace{t_{\alpha_1+\alpha_2},t_{\alpha_1+\alpha_2}}_{\text{2
times }},\notag\\ &\quad
\underbrace{t_{\alpha_1+\alpha_2+1},t_{\alpha_1+\alpha_2+1},t_{\alpha_1+\alpha_2+1}}_{\text{3
times }},\dots).\label{1111}\end{align} Here
$|\alpha|:=\alpha_1+\alpha_2+\dotsm$. Then, for any
$f^{(n)},g^{(n)}\in\mathcal S^{\odot n}_{\mathbb C} $,
 \begin{gather} (f^{(n)},g^{(n)})_{F^{(n)}(L^2({\mathbb R}))}
 =n!\sum_{\alpha\in{\mathbb Z}_+:\, 1\alpha_1+2\alpha_2+\dots=n}\frac{n!}{\alpha_1!\, 1^{\alpha_1}\alpha_2!\,
2^{\alpha_2}\dotsm}\notag\\ \times \int_{{\mathbb
R}^{|\alpha|}}\overline{(D_\alpha
f^{(n)})(t_1,\dots,t_{|\alpha|})}\times(D_\alpha
g^{(n)})(t_1,\dots,t_{|\alpha|}) \, dt_1\dotsm dt_{|\alpha|}.
\label{2222}\end{gather}

The following theorem  \cite{L3,L4} describes the action of
$A_\lambda(\varphi)$ on $\mathcal F_{\mathrm{fin}}(\mathcal S)$.

\begin{theorem}\label{ujftyf}
For each $\lambda\in{\mathbb R}$ and $\varphi\in\mathcal S$, we have
on $\mathcal F_{\mathrm{fin}}(\mathcal S)$: $$
A_\lambda(\varphi)=A^+(\varphi)+\lambda A^0(\varphi)+\mathfrak A^-
(\varphi).$$ Here $\mathfrak A^-(\varphi)$ is the restriction to
$\mathcal F_{\mathrm{fin}}(\mathcal S)$ of the adjoint operator of
$A^+(\varphi)$ in $F(L^2({\mathbb R}))=\bigoplus_{n=0}^\infty
F^{(n)}(L^2({\mathbb R}))n!$, and
$$ \mathfrak A^-(\varphi)=A^-(\varphi)+ A_1^-(\varphi) ,$$
where \begin{equation}\label{jhfhyd} (A^-_1(\varphi)
f^{(n)})(t_1,\dots,t_{n-1})=n(n-1)\big(\varphi(t_1)f^{(n)}
(t_1,t_1,t_2,t_3,\dots,t_n)\big)^\sim.\end{equation}

Furthermore, each $A_\lambda(\varphi)$ is essentially self-adjoint
on $\mathcal F_{\mathrm{fin}}(\mathcal S)$.
\end{theorem}

We see that  $A_\lambda(\varphi)$ has creation, neutral, and
annihilation parts. Therefore, the family of self-adjoint commuting
operators $(A_\lambda(\varphi))_{\varphi\in\mathcal S}$ is a Jacobi
field in $F(L^2({\mathbb R}))$ (compare with  \cite{bere} and
the references therein).

From Theorem \ref{ujftyf}, one concludes that, for any
$f^{(n)}\in\mathcal S_{\mathbb C}^{\odot n}$, the orthogonal
projection ${:}\langle\omega^{\otimes
n},f^{(n)}\rangle{:}_{\lambda}$ has a $\mu_\lambda$-version $\langle
{:}\omega^{\otimes n}{:}_{\lambda},f^{(n)}\rangle$, where
${:}\omega^{\otimes n}{:}_{\lambda}\in\mathcal S'{}^{\odot n}$ are
given by the following recurrence relation:
\begin{gather*}
{:}\omega^{\otimes 0}{:}_{\lambda}=1,\quad  {:}\omega^{\otimes 1}{:}_{\lambda}=\omega,\\
{:}\omega^{\otimes(n+1)}{:}_{\lambda}(t_1,\dots,t_{n+1})=\big(
{:}\omega^{\otimes n}{:}_{\lambda}(t_1,\dots,t_n)\omega(t_{n+1})\big)^\sim \\
\text{}-n \big({:}\omega^{\otimes(n-1)}{:}_{\lambda}(t_1,\dots,t_{n-1})1(t_n)\delta(t_{n+1}-t_{n})\big)^\sim\\
\text{}-\lambda n \big({:}\omega^{\otimes n}{:}_{\lambda}(t_1,\dots,t_{n})\delta(t_{n+1}-t_{n})\big)^\sim\\
\text{}-n(n-1)\big({:}\omega^
{\otimes(n-1)}{:}_\lambda(t_1,\dots,t_{n-1})\delta(t_n-t_{n-1})\delta(t_{n+1}-t_n)\big)^\sim,
\qquad n\in{\mathbb N}.
\end{gather*}

We will now construct a space of test functions. It is possible to
show  \cite{KL} that the Hilbert space $\mathcal F_1(\mathcal S_1)$ is
densely and continuously embedded into $F(L^2({\mathbb R}))$. This
embedding is understood in the sense that $\mathcal F_1(\mathcal
S_1)$  is considered as the closure of $\mathcal
F_{\mathrm{fin}}(\mathcal S)$ in the respective norm. Therefore, the
nuclear space $\mathcal F_1(\mathcal S)$ is densely and continuously
embedded into $F(L^2({\mathbb R}))$. We will denote by $(\mathcal
S)_\lambda ^1$ the image of $\mathcal F_1(\mathcal S)$ under
$I_\lambda$.

Using a result from the theory of test and generalized functions
connected with a generalized Appell system of
polynomials  \cite{KSS,KK},  one proves the following theorem.

\begin{theorem}
For each $\lambda\in{\mathbb R}$,  $$ (\mathcal
S)_{\lambda}^1=\mathcal E_{\mathrm{min}}^1(\mathcal S').$$

\end{theorem}

Thus,  $\mathcal E_{\mathrm{min}}^1(\mathcal S')$ appears to be a
universal space for our purposes.

Taking to notice that $(\mathcal S)_\lambda^1$ is the image of
$\mathcal F_1(\mathcal S)$ under $I_\lambda$, we will identify the
dual space $(\mathcal S')_\lambda^{-1}$ of $(\mathcal S)_\lambda^1$
with $\mathcal F_{-1}(\mathcal S')$. Notice, however, that now the
dual pairing
  between elements of $(\mathcal S')_\lambda^{-1}$ and $(\mathcal S)_\lambda^1$
  is obtained {\it not} through the scalar product in $L^2(\mathcal S',\mu_\lambda)$, or equivalently in $F(L^2({\mathbb R}))$, but through the scalar product in the usual Fock space $\mathcal F(L^2({\mathbb R}))$.  In particular, such a realization of the dual space $(\mathcal S')_\lambda^{-1}$ is convenient for developing Wick calculus on it.

By \eqref{jhfhyd}, the operator $A_1^-(\varphi)$ has the following
representation through the operators $\partial_t$ and
$\partial_t^\dag$:
$$ A_1^-(\varphi)=\int_{\mathbb R}\varphi(t)\partial_t^\dag\partial_t\partial_t\,dt.$$ Therefore, by Theorem \ref{ujftyf}, we get:
$$ A_\lambda(\varphi)=\int_{\mathbb R}\varphi(t)(\partial_t^\dag +\lambda\partial_t^\dag\partial_t+\partial_t+
\partial_t^\dag\partial_t\partial_t)\,dt.$$
Hence, the corresponding quantum white noise, denoted by
$W_\lambda(t)$, is given by
\begin{equation}\label{jfhhyhd} W_\lambda(t)=\partial_t^\dag +\lambda\partial_t^\dag\partial_t+\partial_t+
\partial_t^\dag\partial_t\partial_t.\end{equation}

Realized on the space $(\mathcal S)_\lambda^1$, the operator
$\partial _t$ acts as follows  \cite{L3}:
$$ (\partial_t F)(\omega)=\int_{\mathbb R} \frac{F(\omega+s\delta_t)-F(\omega)}{s}\,\nu_\lambda(ds),\qquad F\in
(\mathcal S)_{\lambda}^1,\ t\in{\mathbb R},\ \omega\in\mathcal S'.
$$

\section{The square of  white noise algebra}
As we have seen in Sections 1 and 2, the Gaussian white noise is
just the sum of the annihilation operator $\partial_t$ and the
creation operator $\partial_t^\dag$, whereas the Poisson white noise
is obtained by adding to the Gaussian white noise a constant times
the product of $\partial_t^\dag$ and $\partial_t$. Let us also
recall that the operators $\partial_t$, $\partial_t^\dag$,
$t\in{\mathbb R}$, satisfy the canonical commutation relations:
\begin{gather} [\partial_t,\partial_s] =
[\partial^\dag_t,\partial^\dag_s]=0 ,\notag\\
[\partial_t,\partial_s^\dag]
=\delta(t-s),\label{z7eawr76}\end{gather} where $[A,B]{:=}AB-BA$.

It was proposed in  \cite{ALV} to develop a stochastic
calculus for higher powers of white noise, in other words, for
higher powers of the operators $\partial_t$, $\partial_t^\dag$. This
problem was, in fact, influenced by the old dream of T.~Hida that
the operators $\partial_t$, $\partial_t^\dag$ should play a
fundamental role in infinite-dimensional analysis.

We will now deal with the squares of $\partial_t$,
$\partial_t^\dag$. The idea is to introduce operators $B_t$ and
$B_t^\dag$ which will be interpreted as $\partial_t^2$ and
$(\partial_t^\dag)^2$, to derive from \eqref{z7eawr76}  the
commutation relations satisfied by $B_t$, $B_t^\dag$, and
$N_t:=B_t^\dag B_t$ and then to consider the quantum white noise
$B_t+B_t^\dag+\lambda N_t$, where $\lambda\in{\mathbb R}$ (compare
with \eqref{jft}, \eqref{jgufgujjj}, and \eqref{gyul}). However,
when doing this, one arrives at the expression $\delta(\cdot)^2$ ---
the square of the delta function. It was proposed in
    \cite{2,3} to carry out a renormalization procedure by
employing the following equality, which may be justified in the
framework of distribution theory:
$$ \delta(\cdot)^2=c\delta(\cdot).$$ Here $c\in{\mathbb C}$ is arbitrary.  This way one gets the following commutation relations:
\begin{gather} [B_t,B_s^\dag]=2c\delta(t-s)
+4\delta(t-s)N_s,\notag\\ [N_t,B_s^\dag]=2\delta(t-s)B_s^\dag,\notag\\
[N_t,B_s]=-2\delta(t-s)B_s,\notag\\
[N_t,N_s]=[B_t,B_s]=[B_t^\dag,B_s^\dag]=0\label{huuh}.
\end{gather}

As usual in mathematical physics,   the rigorous meaning of the
commutation relations \eqref{huuh} is  that they should be
understood in the smeared form. Thus, we introduce the smeared
operators $$B(\varphi):=\int_{{\mathbb R}} \varphi(t) B_t\, dt,\
B^\dag(\varphi):=\int_{{\mathbb R}^d} \varphi(t) B^ \dag_t\,
dt,\quad N(\varphi):=\int_{{\mathbb R}}\varphi(t)N_t\,dt,$$ where
$\varphi\in\mathcal S$, and then the commutation relations between
these operators take the form
\begin{gather} [B(\varphi),B^\dag(\psi)]=2c\langle \varphi,\psi\rangle
+4N(\varphi\psi),\notag\\ [N(\varphi),B^\dag(\psi)]=2B^\dag(\varphi\psi),\notag\\
[N(\varphi),B(\psi)]=-2B(\varphi\psi),\notag\\
[N(\varphi),N(\psi)]=[B(\varphi),B(\psi)]=[B^\dag(\varphi),B^\dag(\psi)]=0,\qquad
\phi,\psi \in \mathcal S . \label{jwaui}\end{gather}
 The operator algebra with generators $B(\varphi),
B^\dag(\varphi),N(\varphi)$, $\varphi\in {\mathcal S}$, and a
central element $1$ with relations \eqref{jwaui} is called the
square of white noise (SWN) algebra.

In    \cite{2}, it was shown that a Fock representation of the
SWN algebra exists if and only if the constant $c$ is real and
strictly positive. In what follows, it will be convenient for us to
choose the constant $c$ to be 2, though this choice is not
essential.

Using the notations of Section \ref{hgyd}, we define
$$ B_t=2(\partial_t+\partial_t^\dag\partial_t\partial_t),\
B_t^\dag =2\partial_t^\dag,\ N_t =2\partial_t^\dag\partial_t. $$
Then it is straightforward to show that the corresponding smeared
operators
\begin{equation} \label{kjggfy}B(\varphi)=2\mathfrak A(\varphi),
\ B^\dag=2A^+(\varphi),\ N(\varphi)=2A^0(\varphi)\end{equation} form
a representation of  a SWN algebra in $F(L^2({\mathbb R}))$. We also
refer to \cite{5} for a unitarily equivalent representation,
see also \cite{SWN}.

Thus, the  quantum white noises $W_\lambda(t)$, $\lambda\in{\mathbb
R}$, (see \eqref{jfhhyhd}) can be thought of as a class  of
(commuting) quantum processes obtained from the SWN algebra
\eqref{kjggfy}.

\end{document}